\documentclass{article}
\usepackage[utf8]{inputenc}
\usepackage{amsmath}
\usepackage{amssymb}
\usepackage{amsthm}
\usepackage{algpseudocode}
\usepackage[hmargin=30mm,vmargin=30mm]{geometry}
\usepackage{graphicx}
\usepackage{hyperref}
\usepackage{mathtools}
\usepackage{todonotes}

\usepackage{color}
\newtheorem{thm}{Theorem}[section]
\newtheorem{lem}[thm]{Lemma}
\newtheorem{rem}[thm]{Remark}
\newtheorem{prop}[thm]{Proposition}
\newtheorem{example}[thm]{Example}

 \newcommand{\DiffCl}{\overline{\operatorname{Diff}}_+}
\newcommand{\Diff}{\operatorname{Diff}_+}
\newcommand{\Imm}{\operatorname{Imm}}

\date{}
\title{Intrinsic Riemannian metrics on spaces of curves: theory and computation.}
\author{Martin Bauer, Nicolas Charon, Eric Klassen, Alice Le Brigant}
\begin{document}
\maketitle

\begin{abstract}
 This chapter reviews some past and recent developments in shape comparison and analysis of curves based on the computation of intrinsic Riemannian metrics on the space of curves modulo shape-preserving transformations. We summarize the general construction and theoretical properties of quotient elastic metrics for Euclidean as well as non-Euclidean curves before considering the special case of the square root velocity metric for which the expression of the resulting distance simplifies through a particular transformation. We then examine different numerical approaches that have been proposed to estimate such distances in practice and in particular to quotient out curve reparametrization in the resulting minimization problems.          
\end{abstract}

\section{Introduction}
Many applications that involve quantitative comparison and statistics over sets of geometric objects like curves often rely on a certain notion of metric on the corresponding shape space. Some of them, such as medical imaging or computer vision, are concerned with the outline of an object, represented by a closed curve, while others, such as trajectory analysis or speech recognition, consider open curves drawing the evolution of a given time process in a certain space, say a manifold. In both cases, it is often interesting when studying these curves to factor out certain transformations (e.g. rotations, translations, reparameterizations), so as to study the shape of the considered object, or to deal with the considered time process regardless of speed or pace.

Beyond computing distances between shapes, a desirable goal in these applications is to perform statistical analysis on a set of shapes, e.g., to compute the mean, perform classification or principal component analysis. For this purpose, considering shapes as elements of a \emph{shape manifold} that we equip with a Riemannian structure provides a convenient framework. In this infinite-dimensional shape manifold, points represent shapes and the distance between two shapes is given by the length of the shortest path linking them -- the geodesic. This approach allows us to do more than simply compute distances: it enables us to define the notion of an optimal deformation between two shapes, and to locally linearize the shape manifold using its tangent space. For instance, given a set of shapes, one can perform methods of standard statistical analyis in the flat representation space given by the tangent space at the barycenter.

The idea of a shape space as a Riemannian manifold was first developed by Kendall \cite{kendall1984shape}, who defines shapes as ``what is left" of a curve after the effects of translation, rotation and changes of scale are filtered out. Mathematically, this means defining the shape space as a quotient space, where the choice of which transformations to quotient out depends on the application. The shapes considered by Kendall are represented by labelled points in Euclidean space and the shape spaces are finite-dimensional. More recent works deal with continuous curves with values in a Euclidean space or a nonlinear manifold, and thus with infinite-dimensional shape spaces. 

There exist two main complementary approaches to define the shape space and its metric. One possibility is to deform shapes by diffeomorphisms of the entire ambient space. In this setting, metrics are defined on the space of spatial deformations, and are called \emph{extrinsic} (or \emph{outer}) metrics as developed in the works of \cite{Grenander1993,trouve98:_diffeom,Beg2005} among other references. Another approach consists in defining metrics directly on the space of curves itself, which are thus called \emph{intrinsic} (or \emph{inner}) metrics. This chapter focuses on the second approach, and studies inner metrics with certain invariance properties. We are specifically interested in the invariance to shape-preserving transformations; in particular to the action of temporal deformations, also called \emph{reparameterizations}, which we represent by diffeomorphisms of the parameter space ($[0,1]$ for open curves, $S^1$ for closed curves). In the following sections, we will introduce a class of invariant Sobolev metrics we call \textit{elastic} on the space of immersed curves which in turn descend to metrics on the space of shapes. These were initially studied in \cite{michor2005vanishing,michor2007overview,Mennucci2008} and in subsequent works. We will then discuss in detail the particular case of the so-called ``Square Root Velocity" (SRV) metric \cite{Jermyn2011}, a first-order invariant metric which allows for particularly simple computations not only for curves in Euclidean spaces but also curves with values in homogeneous spaces or even Riemannian manifolds. Finally, we review different methods to factor out the action of the reparametrization group, which, because of its infinite-dimensionality, presents an important challenge in the computation of distances and geodesics in this framework.

\section{Matching of geometric curves based on reparametrization-invariant Riemannian metrics}

\begin{figure}
    \centering
   \raisebox{.8cm}{ \includegraphics[trim = 50mm 20mm 50mm 15mm ,clip,width=5cm]{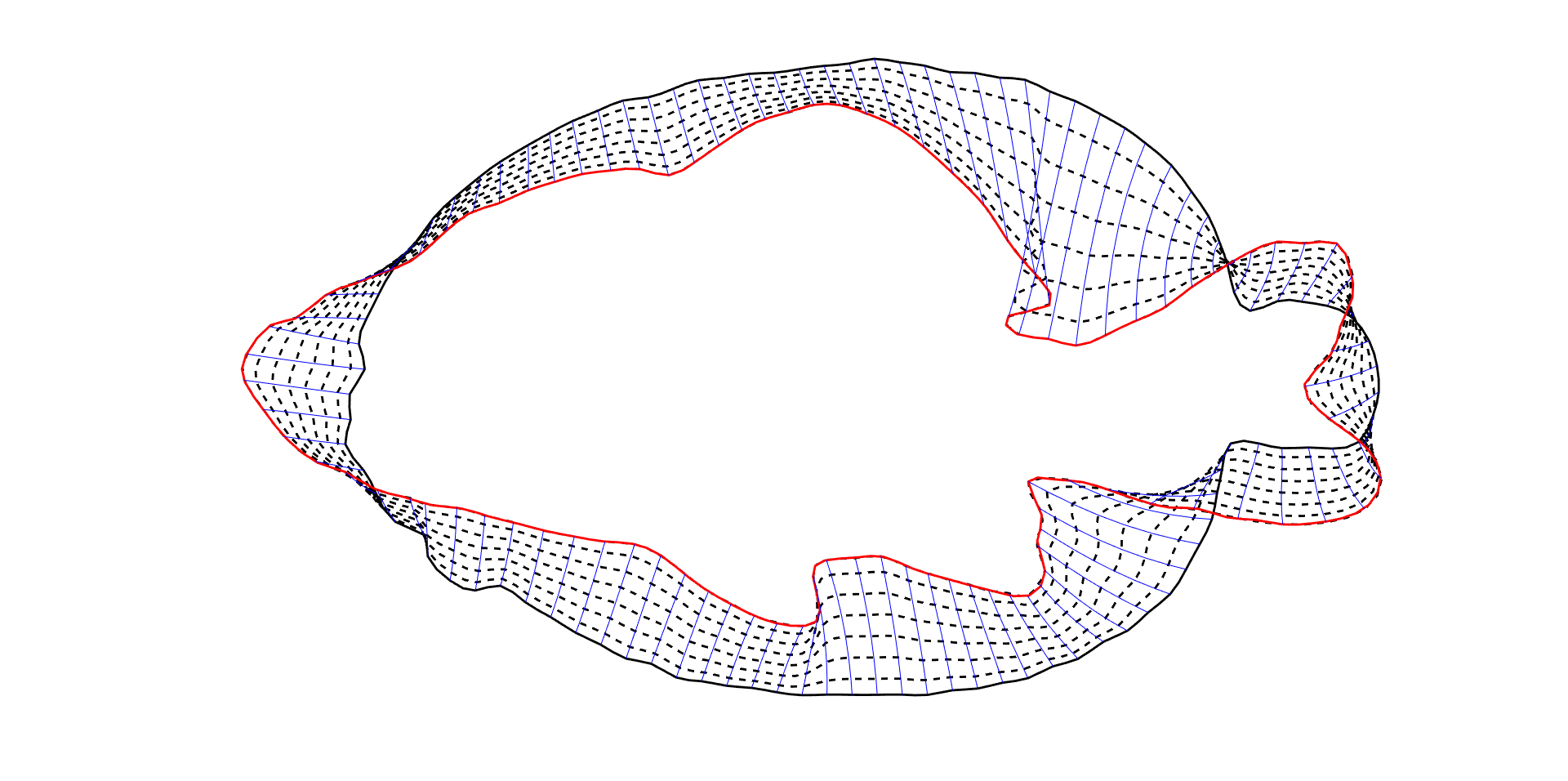}}
    \includegraphics[width=16em]{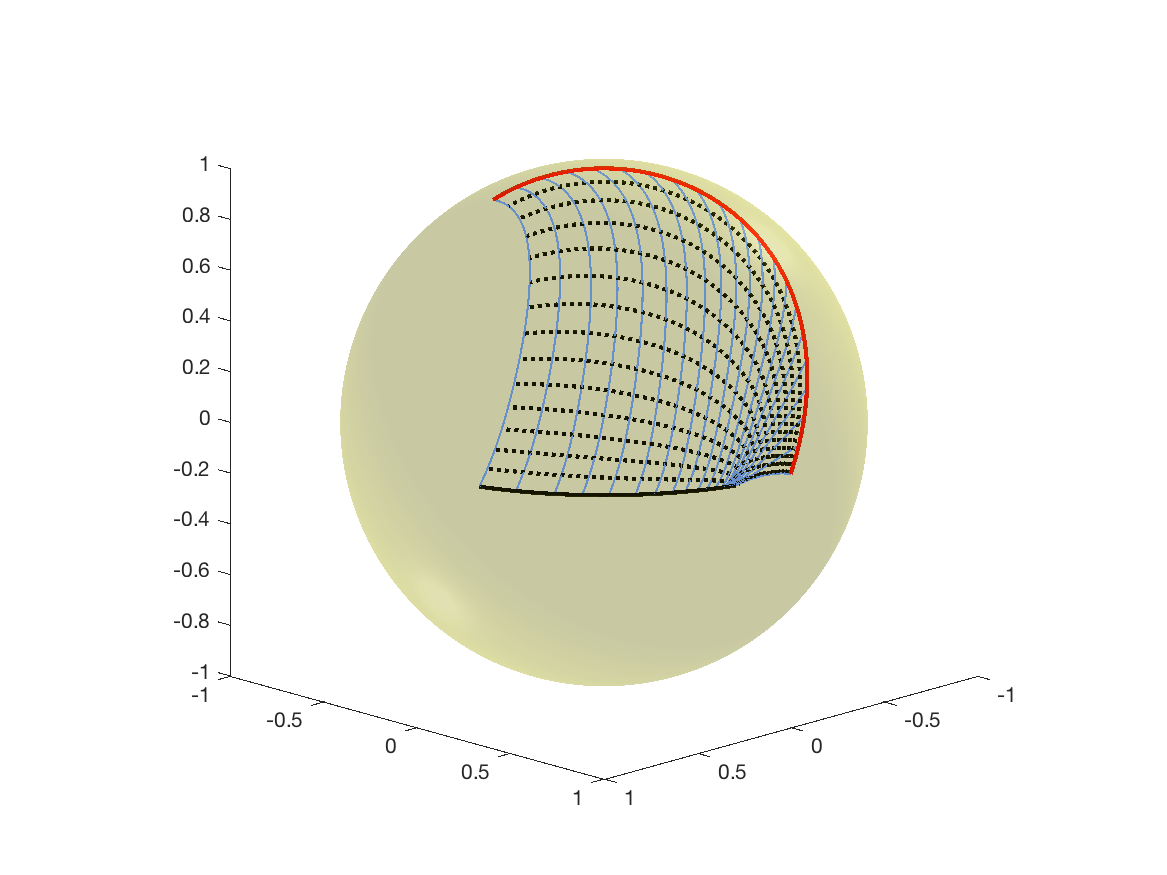}\hspace*{-1em}
   \raisebox{.3cm}{ \includegraphics[width=13em]{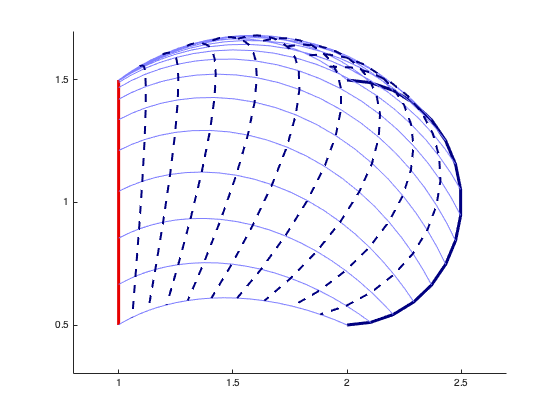}}
    \caption{Examples of geodesics on spaces of unparametrized curves w.r.t elastic metrics  (target curve in red). Some intermediate curves $c(t,\cdot)$ are shown in dashed line and the trajectory of a few specific points in blue. Left figure:  second-order Sobolev metric, estimated with the approach of \cite{bauer2018relaxed}, cf. Section \ref{ssec:relaxed_approach}. Middle figure: SRV-metric for curves with values on homogeneous spaces as implemented in~\cite{su2018comparing}, where the optimal reparametrization is estimated using dynamic programming; cf. Section~\ref{sec:homspaces} and Section~\ref{sec:dynamic}. Right figure: SRV-metric for manifold-valued curves in the hyperbolic plane, as implemented in~\cite{lebrigant2019discrete} with successive horizontalizations; cf. Section~\ref{ssec:SRV_manifold}, method 1 and Section~\ref{sec:horizontalizations}.}
    \label{fig:example_curve_geodesics}
\end{figure}

\subsection{General framework}
Let $D$ be either the interval $I=[0,1]$ or the circle $S^1$ and $(M,\langle.,.\rangle)$ a finite dimensional Riemannian manifold with $TM$ denoting its tangent bundle. In the following we introduce the central object of interest in this book chapter, the infinite dimensional manifold of open (respectively, closed) curves.
\begin{lem}[\cite{michor1980manifolds}]
The space of smooth, regular curves 
\begin{equation}
\operatorname{Imm}(D,M)=\left\{c\in C^{\infty}(D,M): \langle c'(u),c'(u)\rangle_{c(u)}\neq 0,\; \forall u\in D\right\}
\end{equation}
is a smooth Fr\'echet manifold with tangent space at $c$ the set of $C^\infty$ vector fields along $c$, i.e.,
\begin{equation}
T_c\Imm(D,M)=\left\{h\in C^{\infty}(D,TM): h\circ\pi=c\right\}\,,    
\end{equation}
where $\pi: TM\to M$ denotes the foot point projection.
\end{lem}

\begin{figure}
    \centering
    \includegraphics[width=20em]{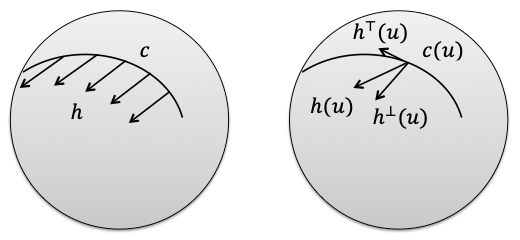}\qquad \vline\qquad
    \raisebox{.7cm}{\includegraphics[width=15em]{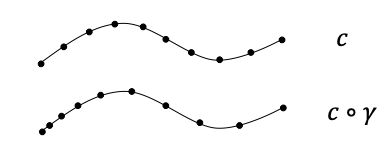}}
    \caption{Left Panel: Tangent vector field to a curve $c(u)$ on the two-dimensional sphere $M=S^2$ (left) and its tangential and normal parts (right). Right Panel: Two different parametrizations of the same geometric curve. }
    \label{fig:schematic_curves}
\end{figure}

The main difficulties for understanding this result stem from the manifold structure of the ambient space $M$. For the convenience of the reader we note that for $M=\mathbb R^d$ the situation simplifies significantly: in that case $\Imm(D,\mathbb R^d)$ is an open subset of the infinite dimensional vector space $C^{\infty}(D,\mathbb R^d)$ and thus tangent vectors to $\Imm(D,\mathbb R^d)$ can be identified with smooth functions with values in $\mathbb R^d$ as well. See Figure \ref{fig:schematic_curves} for a schematic explanation of the involved objects.

In most applications in shape analysis one is not interested in the parametrized curve itself, but only in its features after quotienting out the action of shape-preserving transformations. Therefore, we introduce the reparametrization group of the domain $D$:
\begin{equation}
\Diff(D)=\left\{\gamma\in C^{\infty}(D,D): \gamma \text{ is an orientation preserving diffeomorphism.}\right\}\,.  
\end{equation}
Similarly to the space of immersions, this space carries the structure of an infinite-dimensional manifold. In fact it has even more structure, namely it is an infinite-dimensional Lie group \cite[Section 4]{hamilton1982inverse}. This group acts on the space of immersed curves by composition from the right and this action merely changes the parametrization of the curve but not its actual shape. See Figure \ref{fig:schematic_curves} for an example of different parametrizations of the same geometric curve.

Similarly, we can consider the left action of the group $\operatorname{Isom}(M)$ of isometries  of $M$  on $\Imm(D,M)$. Note that the isometry group is always a finite-dimensional group; e.g. for $M=\mathbb R^d$ the group $\operatorname{Isom}(M)$ is generated by the set of translations and linear isometries\footnote{In some applications one is also interested in moding out the action of the scaling group, which requires a slight modification of the family of elastic metrics. We will not discuss these details here, but refer the interested reader to the literature, e.g.~\cite{Bruveris2017_preprint}.}. Thus, the action of the infinite-dimensional group $\Diff(D)$ is the most difficult to deal with, both from a theoretical and an algorithmic viewpoint. This allows us now to introduce the shape space of curves\footnote{To be mathematically exact, one should limit oneself to the slightly smaller set of free immersions in this definition, as the quotient space has some mild singularities without this restriction. We will, however, ignore this subtlety for the purpose of this book chapter.}
\begin{equation}
\label{eq:def_shape_space}
\mathcal S(D,M):=\Imm(D,M)/\left(\Diff(M)\times \operatorname{Isom}(M)\right))
\end{equation}
 Note: sometimes we use the phrase ``unparametrized shape" to refer to an element of the shape space $\mathcal S(D,M)$ and we shall write $[c] \in \mathcal S(D,M)$ the equivalence class of a parametrized curve $c$.  
\begin{lem}[\cite{cervera1991action,binz1981manifold}]
The shape space $\mathcal S(D,M)$ is a smooth Frechet manifold and the projection $p:\Imm(D,M)\to \mathcal S(D,M)$ is a smooth submersion.
\end{lem}
This means specifically that the mapping $p$ is Frechet-differentiable and that for any $c \in \Imm(D,M)$, $dp(c)$ is onto from $T_c\Imm(D,M)$ to $T_{[c]} \mathcal S(D,M)$. The so called vertical space at $c$ associated to the submersion is defined as $\operatorname{Ver}_c = \{h \in T_c\Imm(D,M) \ | \ dp(c) \cdot h = 0 \}$.  

We aim to introduce Riemannian metrics on the shape space $\mathcal S(D,M)$ by defining metrics on the space of parametrized curves that satisfy certain invariance properties. In the literature these metrics are also referred to as elastic metrics, as they account for both bending and stretching of the curve. 

A Riemannian metric on $\Imm(D,M)$ is a smooth family of inner products $G_c(.,.)$ on each tangent space $T_c\Imm(D,M)$ and we call such a metric $G$ {\em reparametrization-invariant} if it satisfies the relation
\begin{equation}
G_c(h,k)=G_{c\circ\gamma}(h\circ\gamma,k\circ\gamma)    
\end{equation}
for all $c\in \Imm(D,M)$, $h,k\in T_c\Imm(D,M)$ and $\gamma \in \Diff(D)$.

In the following we will introduce the class of Sobolev type metrics. For the convenience of the reader we will first discuss the special case of a first-order metric and $M=\mathbb R^d$. We will then generalize this to the more complicated situation of curves with values in general manifolds and more general metrics. For a curve $c\in \Imm(D,\mathbb R^d)$ and tangent vectors $h,k\in C^{\infty}(D,\mathbb R^d)$ we let
\begin{align}
G_c(h,k)&=\int_D \left(\langle h,k\rangle + \left\langle \frac{h'}{|c'|},\frac{k'}{|c'|}\right\rangle\right)|c'|du
=   \int_D \left(\langle h,k\rangle + \langle D_s h,D_s k\rangle\right)ds\,,
\end{align}
where the desired invariance follows directly by integration using substitution. Here $D_s=\frac{\partial_u}{|c'|}$ and $ds=|c'|du$ denote differentiation and integration with respect to arclength. %
These definitions naturally generalize to curves with values in abstract manifolds by replacing the partial derivative $\partial_u$ in $D_s$ by the covariant derivative with respect to the curve velocity $\nabla_{c'(u)}$. We will denote the induced differential operator as $\nabla_s = \frac{\nabla_{c'}}{|c'|}$. 

Using this notation, a reparametrization-invariant Sobolev metric of order $n$ on the space of manifold valued curves can be defined via
\begin{equation}\label{Sobolevmetric}
G_c(h,k)=\sum_{i=0}^n \int_D \langle \nabla_s^ih, \nabla_s^i k\rangle_{c}\; ds.
\end{equation}
More generally we can consider metrics that are defined by an abstract, positive, pseudo-differential operator $L_c$, that satisfies the equivariance property $L_c(h)\circ \gamma=L_{c\circ \gamma}(h\circ \gamma)$ for all reparametrizations $\gamma$, immersions $c$ and tangent vectors $h$. The corresponding metric can then be written via
\begin{equation}
G_c(h,k)= \int_D \langle L_c(h),
L_c(k) \rangle_{c}\; ds\;.    
\end{equation}
A particularly important example of such metrics is given by the family of elastic $G^{a,b}$-metrics -- first introduced by Mio et.al.~\cite{Mio2007} for the case of planar curves: 
\begin{equation}\label{Gab-metric}
G_c^{a,b}(h,k)=\int_D a^2\langle (\nabla_s h)^{\top},(\nabla_s k)^\top \rangle +b^2
\langle (\nabla_s h)^{\bot},(\nabla_s k)^\bot \rangle ds\,,
\end{equation}
where $a,b>0$ are constants and $\bot$ and $\top$ denote the projection on the normal (respectively, tangential) part of the tangent vector. Here normal and tangential are calculated with respect to the foot point curve $c$, as illustrated in Figure \ref{fig:schematic_curves}.

As a next step we will show that the invariance of the metric $G$ will allow us to define an induced metric on the shape space of unparametrized curves. Before we are able to formulate this result we review some basic facts on  Riemannian submersions. Therefore
let $(\mathcal M,g_1)$ and $(\mathcal{N},g_2)$ be two (possibly infinite dimensional) Riemannian manifolds. A Riemannian submersion is a submersion  $p:(\mathcal M,g_1)\to (\mathcal N,g_2)$ such that in addition $dp: \operatorname{Hor}\to T\mathcal N$ is an isometry. Here $\operatorname{Hor}\subset T\mathcal M$
is the horizontal bundle, which is defined as the $g_1$-orthogonal complement of the vertical bundle
$\operatorname{Ver}:=\operatorname{ker}(dp)\subset T\mathcal M$.
Classical results in Riemannian geometry allow us now to connect the geometry of the two Riemannian manifolds. Most importantly, for our purposes, is the fact that geodesics on $(\mathcal N, g_2)$ correspond to horizontal geodesics on $(\mathcal M,g_1)$. Thus Riemannian submersions are  a convenient construction in our quotient space situation, as it allows, by restricting the calculations to horizontal curves, to perform most of the analysis on the top space, i.e. the space of parametrized curves. 

We are now able to describe the Riemannian submersion picture for the shape space of unparametrized curves.
 Consequently, this gives rise to the following result:
\begin{thm}
The reparametrization-invariant metrics~\eqref{Sobolevmetric}--\eqref{Gab-metric} descend to smooth Riemannian metrics on the quotient space $\mathcal S(D,M)$ such that the projection $p$ becomes a Riemannian submersion. 
\end{thm}
We want to emphasize here that this theorem is non-trivial in our setting:
in finite dimensions, the invariance of the Riemannian metric would always imply the existence of a Riemannian metric on the quotient space, such that the projection is a Riemannian submersion. In our infinite-dimensional situation the proof is slightly more delicate, as one has to show the existence of the horizontal bundle by hand. This can be done by adapting a variant of Moser's trick to the present setting. For the reparamtrization invariant  metrics studied in this chapter, the horizontality condition, requires one essentially to solve a differential equation of order $2n$ with $n$ being the order of the metric. In the case where one is only interested in factoring out the reparametrization group these two subspaces are given by
\begin{align}
&\operatorname{Ver}_c=\left\{h=a.c'\in T_c\Imm(D,M): a\in C^{\infty}(D,\mathbb R)
\right\},\\
&\operatorname{Hor}_c=\left\{k\in T_c\Imm(D,M): G_c(k,ac')=0 \text{ for all } a\in C^{\infty}(D,\mathbb R) 
\right\},
\end{align}
see e.g.~\cite{michor2007overview,Bauer2011b}. If one wants to factor out in addition the group of isometries of $M$, one has to change the definition of the vertical and thus horizontal bundle accordingly. The exact formulas will depend on the manifold $M$.

The above theorem  allows us to develop algorithms on the quotient space $\mathcal S(D,M)$, while performing most of the operations on the space of parametrized curves.
In the following, we discuss how to express the geodesic distance resulting from the above Riemannian metric, that will serve as our similarity measure on the space of shapes. We will first do this for parametrized curves and then in a second step describe the induced distance on the space of geometric curves.
For parametrized curves $c_0,c_1\in\Imm(D,M)$ we have 
\begin{align}\label{distancefunction}
\operatorname{dist}(c_0,c_1)= \operatorname{inf} \int_0^1 \sqrt{G_c(\partial_t c,\partial_t c)} dt\,,   
\end{align}
where the infimum has to be calculated over all paths $c: [0,1]\to\Imm(D,M)$ such that 
$c(0)=c_0$ and $c(1)=c_1$. In the following we will usually view paths of curves as functions of two variables $c(t,u)$ where $t\in[0,1]$ is the time variable along the path and $u\in D$ the curve parameter. 

The induced geodesic distance on the quotient shape space $\mathcal S(D,M)$ can now be calculated via
\begin{align}
\label{distancefunction_quotient}
\operatorname{dist}^{\mathcal S}([c_0],[c_1])=
\underset{\substack{\gamma\in \Diff(D)\\g\in  \operatorname{Isom}(M)}}{\operatorname{inf}} \operatorname{dist}(c_0,g\circ c_1\circ \gamma)=\underset{\substack{\gamma\in \Diff(D)\\g\in  \operatorname{Isom}(M)}}{\operatorname{inf}}  \operatorname{dist}(g\circ c_0\circ\gamma,c_1) \;.
\end{align}
Note that this can be formulated as a joint optimization problem over the path of curves $c$, the reparametrization function $\gamma$ and the isometry $g\in  \operatorname{Isom}(M)$.

In finite dimensions, geodesic distance always gives rise to a true distance function, i.e., it is symmetric, positive and satisfies the triangle inequality. On the contrary, this can fail quite spectacularly in this infinite-dimensional situation, as the geodesic distance can vanish identically on the space. This phenomenon has been found first by Eliashberg and Polterovich for the $W^{-1,p}$-metric on the symplectomorphism group~\cite{eliashberg1993bi}. In the context of reparametrization-invariant metrics on space of immersions this surprising result has been proven by Michor and Mumford~\cite{michor2005vanishing}. In the following theorem we summarize results on the geodesic distance for the class of Sobolev metrics. See \cite{michor2007overview,bauer2012vanishing,bauer2018vanishing,jerrard2019vanishing} and the references therein for further information on this topic.
\begin{thm}
The geodesic distance of the reparametrization-invariant $L^2$-metric -- as defined in equation~\eqref{Sobolevmetric} with $n=0$ -- vanishes on both the space of regular parametrized  curves $\Imm(D,M)$ and on the shape space $\mathcal S(D,M)$. On the other hand, the geodesic distance is positive on both of these spaces if the order of the Sobolev metric is at least one. 
\end{thm}

This result suggests that metrics of order at least one are potentially well-suited for applications in shape analysis. For such applications, one is usually interested in computing numerically the geodesic distance as well as the corresponding optimal path between two given curves. In Riemannian geometry, these optimal paths are called minimizing geodesics and they are locally described by the so-called geodesic equation, which is simply the first-order optimality condition for the length functional as defined in~\eqref{distancefunction}. In our context these equations become rather difficult; they are nonlinear PDEs of order $2n$ (where $n$ is the order of the metric). Nevertheless there exist powerful results on existence of solutions. 

In order to formulate these results we need to introduce the space of all immersions of finite Sobolev regularity, i.e., for $s>\frac{3}{2}$ 
we consider the space 
\begin{equation}
\Imm^s(D,M):=\left\{c\in H^s(D,M): |c'|\neq 0 \right\}\,,   
\end{equation}
which is a smooth Banach manifold. Here $H^s(D,M)$ denotes the Sobolev space of order $s$, see e.g. \cite{bauer2020sobolev} for the exact definition in a similar notation. Note that the condition $|c'|\neq 0$ is well defined as all functions in $H^s(D,\mathbb R^d)$ are $C^1$ for $s>\frac{3}{2}$.  We are now able to state the main result on geodesic and metric completeness, which is of relevance to our applications. In order to keep the presentation as concise as possible we will formulate this result for closed curves and will only comment on the open curves case below. 
\begin{thm}[\cite{bruveris2014geodesic,Bruveris2014b_preprint,bauer2020sobolev}]\label{thm:metric_geod_comp}
Let $\operatorname{dist}$ be the geodesic distance of the Sobolev metric $G$, as defined in ~\eqref{Sobolevmetric}, of order $n\geq 2$ on the space $\Imm(S^1,M)$ of smooth regular curves. The following statements hold:
\begin{enumerate}
    \item The metric $G$ and its corresponding geodesic distance function extend smoothly to the space of Sobolev immersions $\Imm^s(S^1,M)$ for all $s\geq n$.
    \item The space $\Imm^n(S^1,M)$ equipped with the geodesic distance function $\operatorname{dist}$ (of the Sobolev metric of order $n$) is a complete metric space.
    \item For any two curves in the same connected component of $\Imm^n(S^1,M)$ there exists a minimizing geodesic connecting them.
\end{enumerate}
\end{thm}
For open curves it has been shown that the constant coefficient metric as defined in ~\eqref{Sobolevmetric} is in fact not metrically complete~\cite{bauer2018relaxed}. The reason for this is, that one can always shrink down a straight line (open geodesic in the manifold $M$ resp.) to a point using finite energy. One can, however, regain the analogue of the above completeness result for open curves by considering a length weighted version of the Riemannian metric, see \cite{bauer2020sobolev}.

As a direct consequence of the completeness results we  obtain the existence of optimal reparametrizations, i.e., the well-posedness of the matching problem on the space of unparametrized curves. 
To state our main result on existence of optimal reparametrizations we introduce the quotient space of Sobolev immersions modulo Sobolev diffeomorphisms
\begin{equation}
\mathcal S^s(D,M):=\Imm^s(D,M)/\Diff^s(D)/ \operatorname{Isom}(M).  
\end{equation} We have not determined whether this space carries the structure of a manifold. Nevertheless, we can consider the induced geodesic distance on this space and obtain the following completeness  result, which we will formulate again for closed curves only.
\begin{thm}[\cite{Bruveris2014b_preprint}]\label{shapes_space:completeness}
Let $n\geq 2$ and let $\operatorname{dist}$ be the geodesic distance of the Sobolev metric of order $n$ on $\Imm^n(S^1,M)$. Then 
$\mathcal S^n(S^1,M)
$ equipped with the quotient distance $\operatorname{dist}^{\mathcal S}$ is a complete metric space. Furthermore, given two  unparametrized  curves $[c_0],[c_1]\in S^n(S^1,M)$ there exists an optimal reparametrization $\gamma$ and isometry $g$, i.e. the infimum  
\begin{equation}\label{dist_shape}
\operatorname{dist}^{\mathcal S}([c_0],[c_1])=   \underset{\substack{\gamma\in \Diff^n(S^1))\\g\in  \operatorname{Isom}(M)}}{\inf}\operatorname{dist}(c_0,g\circ c_1\circ \gamma)
\end{equation}
is attained. Here $c_0,c_1\in \Imm(S^1,M)$ can be taken as arbitrary representatives of the geometric curves $[c_0]$ and $[c_1]$.
\end{thm}
In the article~\cite{Bruveris2014b_preprint} this result is formulated for the action of the infinite-dimensional group $\Diff(S^1)$ only and for $M=\mathbb R^d$ only. The proof can however be easily adapted to incorporate the action of the compact group $\operatorname{Isom}(M)$ and, using the results of~\cite{bauer2020sobolev}, it directly translates to the case of manifold valued curves. Similar as in theorem~\ref{thm:metric_geod_comp} this results continues to hold for open curves after changing the Riemannian metric to a length weighted version.

For further results on general Sobolev metrics on spaces of curves we refer to the vast literature on the topic, including~\cite{sundaramoorthi2007sobolev,bauer2014overview,klassen2004analysis,michor2007overview,younes1998computable,bauer2019fractional,tumpach2017quotient}. An example of a geodesic between two planar closed curves for a second order Sobolev metric is shown in Figure \ref{fig:example_curve_geodesics} (left), which was computed with the approach described later in Section \ref{ssec:relaxed_approach}. 
In the following section, we will study one particular metric of order one that
will lead to explicit formulas for geodesics and geodesic distance on open, parametrized curves. This will in turn allow us to recover the results on existence of geodesics and optimal reparametrizations. These optimal objects will however fail to have the regularity properties that the optimizers in this section were guaranteed to have.

\subsection{The SRV framework}
\label{ssec:SRV_general}
\subsubsection{Curves in $\mathbb{R}^d$}
\label{ssec:SRV_euclidean}
The reparametrization-invariant Riemannian metrics discussed above are designed to induce Riemannian metrics on the space of shapes. In general, calculating geodesics and distances with respect to these metrics requires numerical optimization, and is often computation-intensive. However, for the case of open curves in $\mathbb R^d$, one of these metrics provides geodesics and distances that are especially easy to compute. This method is known as the ``Square Root Velocity" (SRV) framework.

The main tool in this framework is the map $Q:\Imm(D,\mathbb R^d)\to C^{\infty}(D,\mathbb R^d)$, often referred to in the literature as the SRV transform or function, defined by
\begin{equation}\label{eq:SRVT_Imm}
 Q(c)(u) =  \frac{{c'}(u)}{\sqrt{|{c'}(u)|}}.
\end{equation}
The importance of this map becomes evident in the following theorem by Srivastava et. al.~\cite{Jermyn2011}, which connects it to the $G^{a,b}$-metric \eqref{Gab-metric} for a particular choice of constants $a$ and $b$:
\begin{thm}\label{thm:SRVmetric}
The mapping $Q$ as defined above is an isometric immersion from the space of immersions modulo translations $\operatorname{Imm}(D, \mathbb{R}^d)/\operatorname{Tra}$ with the elastic $G^{1,1/2}$-metric to
$C^{\infty}(D,\mathbb R^d)$ with the flat $L^2$-metric.
\end{thm}
\begin{rem}
\label{rem:simplifying_transforms}
This theorem essentially allows us to transform the computations from a complicated nonlinear manifold to a vector space equipped with a flat metric. In particular, we will see that it leads to explicit formulas for both geodesics and geodesic distance in the case of open curves. 
For planar curves ($d=2$) an analogous transformation for the elastic $G^{a,b}$-metric with $a=b=1$ was found earlier by Younes et. al.~\cite{younes1998computable,Michor2008a}. These transformations have been generalized to all parameters satisfying $a^2-4b^2\geq 0$ (curves in $\mathbb R^d$) by Bauer et. al. in \cite{Bauer2014b} and more recently to arbitrary parameters (planar curves) by Kurtek and Needham~\cite{KuNe2018}. We will
focus in this book chapter solely on the SRV transform, but many of the results are also true for these
other transformations and metrics.
\end{rem}
In the following we will describe the SRV framework in the case of open curves and we will only comment briefly on applications of the SRV transform to closed curves at the end of the section.

\paragraph{Open Curves.}
The reason for treating the case of open curves separately is the fact that the mapping $Q$ becomes a bijection, which will allow us to completely transform all calculations to the image of $Q$ -- a vector space. While we could perform all of these operations in the smooth category, it turns out to be beneficial to consider this method on a much larger space, which will then turn out to be the metric completion of the space of smooth immersions with respect to the SRV-metric.  

Henceforth, for $I=[0,1]$, let $AC(I,\mathbb R^d)$ denote the set of absolutely continuous functions $I\to\mathbb R^d$. Since the considered metric will be invariant under translation, we  standardize all curves to begin at the origin; therefore, let $AC_0(I,\mathbb R^d)$ denote the set of all $c\in AC(I,\mathbb R^d)$ such that $c(0)=0$. We can extend the mapping $Q$ as defined in $\eqref{eq:SRVT_Imm}$ to a mapping on this larger space via $Q:AC_0(I,\mathbb R^d)\to L^2(I,\mathbb R^d)$ as follows
\begin{equation}
 Q(c)(u) = \left\{ \begin{array}{ll}
                \frac{{c'}(u)}{\sqrt{|{c'}(u)|}} & \mbox{if ${c'}(u) \neq 0$};\\
                0 & \mbox{if ${c'}(u) = 0$}.
\end{array} \right. 
\end{equation}
A straightforward calculation shows that $Q$ has an explicit inverse given by
\begin{equation}
c(u)=Q^{-1}(q)(u)=\int_0^u|q(y)|q(y)dy,    
\end{equation}
and, thus, that $Q$ is a bijection. $\Diff(I)$ acts on $AC_0(I,\mathbb R^d)$ from the right by composition; hence, there is a unique right action of $\Diff(I)$ on $L^2(I,\mathbb R^d)$ that makes $Q$ equivariant. The explicit formula for this action is
\begin{equation}\label{L2action}
    (q*\gamma)(u)=\sqrt{\gamma'(u)}q(\gamma(u)),
\end{equation}
where $q\in L^2(I,\mathbb R^d)$ and $\gamma\in \Diff(I)$.
Furthermore, the action of $\Diff(I)$ on $L^2(I,\mathbb R^d)$ defined by \eqref{L2action} is by linear isometries; this follows directly by an application of integration by substitution. Finally, because $Q$ is a bijection, we can use it to induce a Hilbert manifold structure (i.e., a smooth structure and a Riemannian metric) on $AC_0(I,\mathbb R^d)$. Note that this Riemannian metric is exactly the extension of the $G^{1,1/2}$-metric to the space of absolutely continous curves, cf. Theorem~\ref{thm:SRVmetric}.

The central theme of the SRV framework is that the isometry $Q$ enables us to tranform many questions involving the geometry of $AC_0(I,\mathbb R^d)$ to questions involving the well understood geometry of  $L^2(I,\mathbb R^d)$. In particular we obtain the following theorem concerning completeness, geodesics and geodesic distance:
\begin{thm}[\cite{lahiri2015precise,bruveris2016optimal}]\label{GeodesicOpen}
The space of absolutely continous curves equipped with the SRV-metric is a geodesically and metrically complete space. Furthermore, given any curves $c_0,c_1\in AC_0(I,\mathbb R^d)$, the unique minimizing geodesic connecting them is given by
\begin{equation}
c(t,u)=Q^{-1}((1-t)Q(c_0)(u)+t Q(c_1)(u)),
\end{equation}
and thus the geodesic distance between $c_0$ and $c_1$ can be calculated via
\begin{equation}
\operatorname{dist}(c_0,c_1)= \sqrt{\int_0^1 |Q(c_0)(u)- Q(c_1)(u)|^2 du }\;.   
\end{equation}
\end{thm}

\paragraph{Optimal Reparametrizations.} At this point, it remains to discuss the existence of optimal matchings in the definition of the quotient metric, namely, given two curves $c_0$, $c_1$, does there exist a reparametrization $\gamma \in \Diff(I)$ that attains the infimum in~\eqref{distancefunction_quotient}? The first result in this direction was obtained by Trouv\'{e} and Younes in \cite{trouve2000class} (we also refer to the discussion in \cite[Section 12.7.4]{younes2019shapes}). In this work they analyze the existence of minimizers for a general class of optimization problems on the group of diffeomorphisms of $[0,1]$.
In the case of the elastic $G^{a,b}$-metrics~\eqref{Gab-metric} for open planar curves and when $a>b$, it implies that the existence of an optimal reparametization $\gamma$ always holds for piece-wise $C^1$ curves. When $a=b$ one needs to assume in addition that there does not exist a flat region of one curve together with a point on the other curve for which the tangent vectors are pointing in opposite directions (and with parameters within a certain distance of one another). However, for $a<b$, the conditions become much more restrictive, as one needs to exclude the situation in which there is an open interval in the parameter domain of one curve where the angle between the tangents and the tangent at a point of nearby parameter in the other curve exceeds $a\pi/b$. In particular, for the SRV-metric, this basically constrains angles between tangent vectors of the two curves to be smaller than $\pi/2$, which is an impractical assumption in typical applications. As we discuss next, it turns out that by allowing instead of a single diffeomorphism a pair of ''generalized" reparametrization functions one can recover an existence result for fairly general classes of curves.

In the following we aim to describe this construction, which will require us to consider the closure of the  $\Diff(I)$-orbits on $AC_0(I,\mathbb R^d)$.
Hence, we define an equivalence relation on $AC_0(I,\mathbb{R}^d)$ by $c_1\sim c_2$ if and only if the $\Diff(I)$ orbits of $Q(c_1)$ and $Q(c_2)$ have the same closure in $L^2(I,\mathbb{R}^d)$. We then define the {\em shape space} of open curves in $\mathbb{R}^d$ as 
$$\mathcal{S}(I,\mathbb{R}^d)=AC_0(I,\mathbb{R}^d)/\sim,$$
and for $c\in AC_0(I,\mathbb{R}^d)$ we let $[c]$ denote the equivalence class of $c$ under $\sim.$

In order to better understand these equivalence classes, we need an expanded version of $\Diff(I)$. To be precise, define 
$\DiffCl(I)$ to be the set of all absolutely continuous functions $\gamma:I\to I$ such that $\gamma(0)=0$, $\gamma(1)=1$, and $\gamma'(u)\geq 0$ almost everywhere. Note that $\DiffCl(I)$ is only a monoid, not a group, since the only elements of $\DiffCl(I)$ that have inverses are those $\gamma$ such that $\gamma'(u)\neq 0$ almost everywhere. We then have the following description of a general equivalence class of $AC_0(I,\mathbb{R}^d)$ under the relation $\sim$:
\begin{lem}[\cite{lahiri2015precise}]\label{EquivChar}
Let $c\in AC_0(I,\mathbb{R}^d)$, and assume that $c'(u)\neq 0$ almost everywhere. Then the equivalence class of $c$ under $\sim$ is equal to
$$\{c\circ\gamma:\gamma\in\DiffCl(I)\}.$$
\end{lem}
Note that if $c'(u)=0$ on a set of nonzero measure, then we cannot directly use Lemma \ref{EquivChar} to characterize $[c]$; however, we can reparametrize $c$  by arclength to obtain another element $\tilde c$ in the same equivalence class as $c$, and then use Lemma \ref{EquivChar} to characterize $[c]=[\tilde c]$.  

We can now define a distance function on the shape space as follows: if $[c_1]$ and $[c_2]$ are elements of $\mathcal{S}(I,\mathbb{R}^d)$, then we let
$$\operatorname{dist}^{\mathcal S}([c_0],[c_1])=\inf_{w_0\in[c_0],w_1\in[c_1]} \|Q(w_0)-Q(w_1)\|_{L^2}.$$
Note that it seems at first that we need to consider reparametrizations of both $c_0$ and $c_1$, because $\DiffCl(I)$ is not a group but only a monoid. However, it can be shown that the infimum will be the same if we only consider reparametrizations of one of the curves. See \cite{lahiri2015precise,bruveris2016optimal}.
The optimal reparametrization problem for curves in $AC_0(I,\mathbb{R}^d)$ can now be formulated as follows: suppose $c_0$ and $c_1$ are elements of $AC_0(I,\mathbb{R}^d)$, and that both have non-vanishing derivatives almost everywhere. Do there exist $\gamma_0$ and $\gamma_1$ in $\DiffCl(I)$ such that
$$\|Q(c_0\circ\gamma_0)-Q(c_1\circ\gamma_1)\|_{L^2}=\operatorname{dist}^{\mathcal S}([c_0],[c_1])\;?$$
The following theorem  gives the known results about this problem.
\begin{thm}[\cite{lahiri2015precise,bruveris2016optimal}]\label{OptReparam}
Let $c_0$ and $c_1$ be elements of $AC_0(I,\mathbb{R}^d)$ with both having non-vanishing derivatives almost everywhere. We have:
\begin{enumerate}
    \item if at least one of these curves is piecewise linear, then a pair $\gamma_0, \gamma_1$ of optimal reparametrizations exists;
    \item if $c_0$ and $c_1$ are both of class $C^1$, then a pair $\gamma_0, \gamma_1$ of optimal reparametrizations exists;
    \item there exists a pair $c_0,c_1\in AC_0(I,\mathbb{R}^d)$, both Lipschitz, for which no pair of optimal reparametrizations exists.
\end{enumerate}
\end{thm}

\begin{rem}\label{rem:explicit_optimalreparams}
Later in this chapter numerical techniques for approximating optimal reparametrizations are discussed.  However, we note here that in \cite{lahiri2015precise} an algorithm is developed for determining precise optimal reparametrizations for the case in which both $c_0$ and $c_1$ are piecewise linear curves. Nevertheless, since this algorithm is computationally rather expensive, usually the numerical methods described in Section~\ref{sec:numerics} are used to solve the matching problem in practice. Furthermore, all of the algorithms that we discuss in Section~\ref{sec:numerics}
solve only for one reparametrization function (as opposed to a pair of optimal reparametrization functions as required by the above theorem). Thus the existence of minimizers for these algorithms is only guaranteed for metrics of order two or higher (by the results of Theorem~\ref{shapes_space:completeness}). For lower order metrics, such as the SRV-metric, the computed distances can approximate the true geodesic distances of arbitrary precision by the density of $\Diff(I)$ in $\DiffCl(I)$.
\end{rem}

\paragraph{Closed curves.}
For applications in which curves correspond to boundaries of planar regions, the SRV framework can be adapted to the space of closed curves. A priori, it is natural to describe a closed curve as an immersion of the circle $S^1$ into $\mathbb{R}^d$; then the natural group of reparametrizations is $\Diff(S^1)$. However, in order to apply the SRV methods already outlined, we will work again in the absolutely continuous category and describe a closed curve by an open curve whose initial and endpoints happen to coincide. 
Hence, we define the set of absolutely continuous, closed curves by  $$AC_0(I,\mathbb{R}^d)_{cl}=\{c\in AC_0(I,\mathbb{R}^d) : c(0)=c(1)\},$$
which is a codimension $d$ submanifold of $AC_0(I,\mathbb{R}^d)$.
In order to endow $AC_0(I,\mathbb{R}^d)_{cl}$ with a Riemannian structure, we simply restrict the SRV-metric on $AC_0(I,\mathbb{R}^d)$ to this submanifold. Unfortunately, $AC_0(I,\mathbb{R}^d)_{cl}$  is not a geodesically convex submanifold, so computing geodesics and geodesic distances is not as straightforward as it is in $AC_0(I,\mathbb{R}^d)$. 

Fortunately, the necessary analytical tools have been developed to solve this problem. To find a geodesic between two curves $c_0$ and $c_1$ in $AC_0(I,\mathbb{R}^d)_{cl}$, on can use the following procedure:
\begin{enumerate}
    \item Calculate a geodesic $\{c_t\}$ between $c_0$ and $c_1$ in $AC_0(I,\mathbb{R}^d)$ using Theorem \ref{GeodesicOpen}.
    \item For each $t\in[0,1]$, project $c_t$ to a nearby point $\tilde c_t$ in $AC_0(I,\mathbb{R}^d)_{cl}$. This requires a gradient algorithm as described in \cite{Jermyn2011,Srivastava2016}.
    \item Deform $\{\tilde c_t\}$ to a geodesic in $AC_0(I,\mathbb{R}^d)_{cl}$ using a path-straightening procedure, as described in \cite{Jermyn2011,Srivastava2016}.
\end{enumerate}
In practice, Step 3 is often omitted to save computation, because the path produced by Step 2 is generally very close to a geodesic. 
In order to find optimal reparametrizations for a pair of closed curves, it is not enough to consider the methods developed for open curves, because of the freedom to choose any point on a closed curve to be its starting and ending point (i.e. the point $c(0)=c(1)$). To remedy this, the algorithms discussed for open curves need to be implemented along a densely spaced set of points on one of the curves in order to choose the matching that leads to the shortest geodesic between the curves. For details, see \cite{Srivastava2016}.

\subsubsection{Curves in Lie groups}\label{sec:homspaces}
In the following sections we will discuss methods for extending the SRV framework to curves in Lie groups, homogeneous spaces and manifolds. We start by the simplest generalization:  curves with values in Lie groups, for which the existence of a designated tangent space, the Lie algebra, makes the generalization of the SRV framework straightforward, cf. \cite{celledoni2016shape,su2018comparing}.

Consider a finite-dimensional Lie group ${\mathfrak G}$ with Lie algebra ${\mathfrak g}=T_e{\mathfrak G}$, where $e\in {\mathfrak G}$ denotes the neutral element. We will assume that ${\mathfrak g}$ has been equipped with an inner product and that this inner product has been extended to a left-invariant Riemannian metric on ${\mathfrak G}$. Following the square root velocity framework (SRVF) described above for curves in $\mathbb R^d$, we define the map 
\begin{align}\label{eq.Q.map}
\begin{cases}
Q: AC(I, \mathfrak{G})\to \mathfrak{G}\times L^2(I,\mathfrak{g})\\
Q(c)=(c(0), q),
\end{cases}
\end{align}
where 
\begin{align}\label{eq.q.map}
q(u)=\left\{  \begin{array}{lcr}
\dfrac{dL_{c(u)^{-1}}c'(u)}{\sqrt{\|c'(u)\|}} &c'(u)\neq 0\\
0       &c'(u)=0
\end{array} \right.
\end{align}

Note that $L_{c(u)^{-1}}$ denotes left translation on ${\mathfrak G}$ by $c(u)^{-1}$, which is added to transport the whole curve to the same tangent space ${\mathfrak g}$. Note also that the second part of this transformation is simply the generalization of the SRV transform for curves in a Euclidean space to curves with values in a Lie group and the first factor is added to keep track of the starting point. In \cite{su2018comparing}, it is shown that the map $Q$ is a bijection.

We put a product metric on ${\mathfrak G}\times L^2(I,{\mathfrak g})$ coming from the left-invariant metric on ${\mathfrak G}$ and the $L^2$-metric on $L^2(I,{\mathfrak g})$. Then the smooth structure and Riemannian metric on ${\mathfrak G}\times L^2(I,{\mathfrak g})$ is pulled back to $AC(I,{\mathfrak G})$ leading to the following explicit formula for the corresponding geodesic distance 
\begin{equation}
\label{eq:dist_SRVF_Lie_groups}
\operatorname{dist}(c_0,c_1)^2=\operatorname{dist}^{\mathfrak G}(c_0(0),c_1(0))^2+\int_0^1 \|q_1(u) -q_0(u)\|^2 du,
\end{equation}
with $\operatorname{dist}^{\mathfrak G}$ being the geodesic distance on the finite dimensional group $\mathfrak G$ and $q_i(u)$ being the $q$-map, as defined in equation~\eqref{eq.q.map}, of the curve $c_i$.
Note, that the smooth structure and metric are invariant under the action of $\Diff(I)$ and also under the left action of ${\mathfrak G}$. For the relation of the corresponding Riemannian metric to the class of elastic metrics as defined in equation~\eqref{Sobolevmetric} we refer to the articles ~\cite{su2018comparing,celledoni2016shape}. 
\begin{example}
To make the above more explicit on a simple example, consider $\mathfrak G$ the Lie group $\operatorname{SO}(n,\mathbb R)$ of real $n\times n$ orthogonal matrices with determinant one, the group operation being the standard matrix product. On the corresponding Lie algebra, which is the space of antisymmetric $n \times n$ matrices, we consider the inner product 
\begin{equation}
\langle A, B\rangle =\operatorname{tr}(A^TB) = -\operatorname{tr}(AB).    
\end{equation}
For a curve $c$ in $\operatorname{SO}(n,\mathbb R)$, its q-map then simply writes:
\begin{equation}
q(u) = \frac{c(u)^{-1}c'(u)}{\sqrt{\operatorname{tr}\Big((c(u)^{-1}c'(u))^T(c(u)^{-1}c'(u)\Big)}} = \frac{c(u)^{T}c'(u)}{\sqrt{\operatorname{tr}(c'(u)^Tc'(u))}}\;.    
\end{equation}
For the last equality we used that $c(u)^{T}=c(u)^{-1}$.
Moreover the geodesic distance on $\mathfrak G$ is given explicitly by $\operatorname{dist}^{\mathfrak G}(c_0,c_1) = \|\log(c_0^{T}c_1)\|_F^2$, where $\log$ denotes the standard matrix logarithm and $\|.\|_F$ the Frobenius norm. This leads to the following specific expression of the SRV distance \eqref{eq:dist_SRVF_Lie_groups} for parametrized curves in $\operatorname{SO}(n,\mathbb R)$:
\begin{equation}
\operatorname{dist}(c_0,c_1)^2=\|\log(c_0(0)^{T}c_1(0))\|_F^2+\int_0^1 \|q_1(u) -q_0(u)\|_F^2 du .   
\end{equation}
\end{example}

\subsubsection{Curves in homogenous spaces}\label{sec:homogeneous}

For homogenous spaces the situation becomes slightly more complicated and will require an additional minimization over a finite-dimensional group. 

We first recall the definition of a homogenous space. 
A {\em homogeneous space} $M={\mathfrak G}/{\mathfrak K}$ is a quotient of a Lie group ${\mathfrak G}$ by a closed Lie subgroup ${\mathfrak K}$. Note that this quotient is interpreted only as a set of left cosets; it cannot be thought of as a quotient group, since there is no assumption that ${\mathfrak K}$ is a normal subgroup. For purposes of this chapter, we will assume that the subgroup ${\mathfrak K}$ is compact. Examples of homogeneous spaces include spheres, Grassmannians, hyperbolic spaces, and spaces of symmetric positive definite matrices which occur in many applications. 

There is a natural left action of ${\mathfrak G}$ on $M={\mathfrak G}/{\mathfrak K}$ and we endow $M$ with a Riemannian metric that is invariant under this ${\mathfrak G}$-action as follows. First, we put a Riemannian metric on ${\mathfrak G}$ that is left-invariant under the action of ${\mathfrak G}$ and bi-invariant under the action of ${\mathfrak K}$. This is always possible using an averaging argument and the compactness of ${\mathfrak K}$. This metric then descends to a metric on $M$ that is invariant under the left action of ${\mathfrak G}$. In
order to study the shape space of curves with values in the homogeneous space $M$, we wish to put a Riemannian metric on the space $AC(I,M)$ that is invariant under the action of $\Diff(I)$ and the natural left action of ${\mathfrak G}$. We now summarize how this is accomplished using a natural adaptation of the SRV approach for Lie-groups from the previous section; see 
\cite{celledoni2016bshape,su2018comparing} for more details. 

The main idea  is to lift curves in $M$ to curves in ${\mathfrak G}$ that are horizontal (i.e., orthogonal to each $\mathfrak{K}$-coset that they meet). This allows us then to use the ideas for curves in Lie groups, as described in the previous section. Therefore let ${\mathfrak k}\subset {\mathfrak g}$ to be the Lie algebra of ${\mathfrak K}$, and let ${\mathfrak k}^\perp$ be the orthogonal complement of ${\mathfrak k}$ in ${\mathfrak g}$. Let $\pi:
\mathfrak{G}\to M$ denote the natural surjection. If we restrict $Q^{-1}$ (the inverse of the map defined in equation (\ref{eq.Q.map})) to ${\mathfrak G}\times L^2(I,{\mathfrak k}^\perp)$, and then compose with $\pi$, we obtain a surjection
$${\mathfrak G}\times L^2(I,{\mathfrak k}^\perp)\to AC(I,M).$$ 
This surjection is not a bijection, because a curve $c$ in $AC(I,M)$ does not have a unique horizontal lift to $\mathfrak{G}$. Rather, it has a unique horizontal lift starting at each point of $\pi^{-1}(c(0))$. To fix this, we define a right action of $\mathfrak{K}$ on ${\mathfrak G}\times L^2(I,{\mathfrak k}^\perp)$ by
$$(c_0,q)*y=(c_0y,y^{-1}qy),$$
where $y\in\mathfrak{K}$, $c_0\in\mathfrak{G}$, and $q\in L^2(I,{\mathfrak k}^\perp)$. Taking the quotient under this action precisely remedies the lack of injectivity, yielding a bijection
$$({\mathfrak G}\times L^2(I,{\mathfrak k}^\perp))/\mathfrak{K}\to AC(I,M).$$
This bijection, which is equivariant with respect to the left action of $\mathfrak{G}$, is the key tool that we use to define a Riemannian metric on $AC(I,M)$. To see this, note first that we can endow ${\mathfrak G}\times L^2(I,{\mathfrak k}^\perp)$ with the natural product metric in the same way that we did in the case of Lie groups. Then, note that this metric is invariant under the right action of $\mathfrak{K}$, so it induces a metric on the quotient space $({\mathfrak G}\times L^2(I,{\mathfrak k}^\perp))/\mathfrak{K}$ and, hence, on $AC(I,M)$. Furthermore, this Riemannian metric is invariant under the left action of $\mathfrak{G}$.

\noindent {\bf Geodesics:} Geodesics in $L^2(I,{\mathfrak k}^\perp)$ are simply straight lines. Let us assume that we can compute geodesics in $\mathfrak{G}$, as well. Then geodesics in $\mathfrak{G}\times L^2(I,{\mathfrak k}^\perp)$ are products of geodesics in these two spaces. To compute geodesics and geodesic distance in $AC(I,M)$, we need to compute geodesics in 
$(\mathfrak{G}\times L^2(I,{\mathfrak k}^\perp))/\mathfrak{K}$. This is accomplished as follows. Suppose we are given two elements of $(\mathfrak{G}\times L^2(I,{\mathfrak k}^\perp))/\mathfrak{K}$, $[(c_1,q_1)]$ and $[(c_2,q_2)]$. In order to calculate a geodesic between them, we must find $y\in\mathfrak{K}$ that minimizes $d((c_1,q_1),(c_2y,y^{-1}q_2y))$. Note that this is a minimization problem over the compact Lie group $\mathfrak{K}$. In fact, the gradient of this function on $\mathfrak{K}$ can be explicitly calculated (see Lemma 5 of \cite{su2018comparing} for the computation), reducing the computation of geodesics to an optimization problem on a compact Lie group with an explicit gradient. This technique yields efficiently computable formulas for geodesics and geodesic distances, see~\cite{celledoni2016bshape,su2018comparing}. See Figure~\ref{fig:example_curve_geodesics_3D} for an example of geodesics between curves on the sphere. Furthermore, analogues of the optimal reparametrization results, cf. Theorem~\ref{OptReparam}, have been proven, see~\cite{su2018comparing}.

Finally, we note that under the framework just described, the Lie group $\mathfrak{G}$ acts on $AC(I,M)$ by isometries. Hence, for some applications, one may wish to mod out by this action (in addition to the reparametrization group) when defining the shape space of open curves in $M$. We observe that the current framework extends very naturally to performing the additional optimization implied by this quotient operation. We refer the reader to \cite{su2018comparing} and \cite{SBK2017} for more details.

\subsubsection{Curves in Riemannian manifolds}\label{ssec:SRV_manifold}

Let us focus again on open curves, i.e. when $D$ is the interval $I=[0,1]$. For manifold-valued curves, the generalization of the SRV framework is no longer straightforward. Here we discuss three different generalizations. The first method builds on the elastic $G^{1,1/2}$-metric, replacing ordinary derivatives by covariant derivatives with respect to the connection $\nabla$ of the base manifold $M$. The two other methods, while not implementing the precise elastic method, are less computationally expensive, and often yield useful comparisons between curves. Both of these methods replace each curve in the Riemannian manifold $M$ by a curve in a single tangent space of $M$, thus moving the computations to that tangent space, while in the first one, computations are done directly in the base manifold.

\paragraph{Method 1.} In the case of curves with values in a general manifold, the elastic $G^{1,1/2}$-metric is no longer a flat metric. However it can still be obtained as a pullback by the SRV transform of a natural metric on the tangent bundle $T\Imm(I,M)$, namely a pointwise version of the Sasaki metric on $TM$. Recall that the Sasaki metric is a natural choice of metric on the tangent bundle $TM$ that depends on the horizontal and vertical projections of each tangent vector. Intuitively, the horizontal projection of a tangent vector of $T_{(p,w)}TM$ for any $(p,w)\in TM$ corresponds to the way it moves the base point $p$, and its vertical projection, to the way it linearly moves $w$. More precisely, define just as in the Euclidean case the SRV transform to be $Q:\Imm(I,M)\rightarrow T\Imm(I,M)$,
\begin{equation*}
Q(c)(u)=c'(u)/\sqrt{|c'(u)|}. 
\end{equation*}
Consider the following metric on the tangent bundle $T\Imm(I,M)$: for any pair $(c,h)\in T\Imm(I,M)$, and any infinitesimal deformations $\xi_1, \xi_2 \in T_{(c,h)}T\Imm(I,M)$ of the pair $(c,h)$, define
\begin{equation}
\label{eq:sasakimetric}
\hat G_{(c,h)}\left(\xi_1,\xi_2\right) = \langle\xi_1(0)^{\operatorname{hor}}, \xi_2(0)^{\operatorname{hor}}\rangle + \int_I \langle\xi_1(u)^{\operatorname{ver}}, \xi_2(u)^{\operatorname{ver}}\rangle \mathrm du,
\end{equation}
where $\xi_1(u)^{\operatorname{hor}}\in TM$ and $\xi_1(u)^{\operatorname{ver}}\in TM$ are the horizontal and vertical projections of the tangent vector $\xi_1(u)\in T_{(c(u),h(u))}TM$ for all $u\in I$. Then, the elastic $G^{1,1/2}$-metric is the pullback of $\hat G$ with respect to the SRV transform $Q$, i.e.
\begin{align}
\label{eq:elasticmanifold}
G^{1,1/2}_c(h,k) = \hat G_{Q(c)}\left( T_cQ(h) , T_cQ(k) \right)
= \langle h(0),k(0)\rangle+\int_I \langle \nabla_{h(u)}Q(c),\nabla_{k(u)}Q(c)\rangle \mathrm du,
\end{align}
for any curve $c\in \Imm(I,M)$ and $h,k \in T_c\Imm(I,M)$, where $\nabla_{h(u)}Q(c)$ denotes the covariant derivative in $M$ of the vector field $Q(c)$ in the direction of the vector field $h$. Notice that here we add a position term to the integral definition \eqref{Gab-metric} of the $G^{1,1/2}$-metric in order to take into account translations. Accordingly, the energy of a path of curves $[0,1]\ni t \mapsto c(t)$ which SRV transform we write $q(t,\cdot) = Q(c(t))$ for the $G^{1,1/2}$-metric is given by
\begin{equation}
\label{eq:energy1}
E(c) = \int_0^1\left(|\partial_tc(t,0)|^2+\int_I |\nabla_tq(t,u)|^2\mathrm du\right)\mathrm dt.
\end{equation}
Here, $\nabla_th$ denotes the covariant derivative in $M$ of a vector field $t\mapsto h(t,u)$ along a curve $t\mapsto c(t,u)$, i.e. $\nabla_th=\nabla_{\partial_tc}h$. A variational approach yields the following conditions for such a path to be geodesic.
\begin{prop}[\cite{lebrigant2016computing}]
\label{prop:geod_equations_Riemannian}
A path of curves $[0,1]\ni t \mapsto c(t)$ is a geodesic for the $G^{1,1/2}$-metric if and only if its SRV representation $q(t)=Q(c(t))$ verifies the following equations:
\begin{align*}
\nabla_t\partial_tc(t,0) + r(t,0)  =& 0 , \quad \forall t\in[0,1], \\
\nabla_t^2 q(t,u) + |q(t,u)| \left( r(t,u) + r(t,u)^T \right) =& 0 , \quad \forall (t,u)\in[0,1]\times I,
\end{align*}
where the vector field $r$ depends on the curvature tensor $\mathcal R$ of the base manifold $M$ and on the parallel transport $\partial_tc(t,v)^{v,u}$ of the vector field $\partial_tc(t,\cdot)$ along $c(t,\cdot)$ from $c(t,v)$ to $c(t,u)$:
\begin{equation*}
r(t,u) = \int_u^1 \mathcal R(q,\nabla_tq)\partial_tc(t,v)^{v,u} \mathrm dv.
\end{equation*} 
\end{prop}
In the flat case $M=\mathbb R^d$, the curvature term $r$ in the geodesic equation vanishes and we obtain $\nabla_t\partial_tc(t,0)=\partial_t^2c(t,0)=0$, $\nabla_t^2q(t,u)=\partial_t^2q(t,u)=0$ for all $(t,u)\in [0,1]\times I$. We then recover the fact that the geodesic for the SRV-metric between two curves in $\mathbb R^d$ links their starting points with a straight line and linearly interpolates between their SRV representations. In the general case, the initial value problem for geodesics can be solved by finite differences, and the boundary value problem by geodesic shooting. In the case where the base manifold $M$ has constant sectional curvature, e.g. the sphere or the hyperbolic plane, a comprehensive discrete framework was proposed in \cite{lebrigant2019discrete} that correctly approximates the continuous setting and makes numerical computations easier.

\paragraph{Method 2.} An important complication linked to curves taking their values in a non-linear manifold is that tangent vectors $h\in T_c\Imm(I,M)$, which are smooth vector fields along the curve $c$, are functions taking their values in different linear spaces. In order to bypass this difficulty, another way to go is to parallel transport the SRV transform of each curve to a single tangent space, namely the tangent space to the curve's starting point. Consider the vector bundle $\pi:{\mathcal C}\to M$ in which the fiber over each point $x\in M$ is the set of smooth functions in the tangent space $T_xM$, i.e. $\pi^{-1}(x)=C^\infty(I,T_xM)$. Then, define a map
\begin{equation*}
\begin{cases}
Q^\parallel:\Imm(I,M)\to{\mathcal C}\\
Q^\parallel(c)=q^\parallel\in\pi^{-1}(c(0)),
\end{cases}
\end{equation*}
where, for each $u\in I$, $q^\parallel(u)$ is obtained by parallel translating the vector ${c'(u)}/\sqrt{|c'(u)|}$ along the curve $c$ from $c(u)$ to $c(0)$. The function $q^\parallel=Q^\parallel(c)$ therefore takes its values in $T_{c(0)}M$ and is called the "Transported Square Root Velocity" (TSRV) representation of the curve $c$. The vector bundle $\mathcal{C}$ is endowed with a metric that, just like \eqref{eq:sasakimetric}, is a pointwise version of the Sasaki metric, i.e. defined for each $(x,v)\in\mathcal C$ and tangent vectors $(w_1,\eta_1),(w_2,\eta_2)\in T_{(x,v)}\mathcal C$ by%
\begin{equation*}
    \hat G_{(x,v)}((w_1,\eta_1),(w_2,\eta_2))=\langle w_1, w_2\rangle + \int_I\langle \eta_1(u),\eta_2(u)\rangle \mathrm du.
\end{equation*}
It is easily shown that the pullback of this metric to $\operatorname{Imm}(I,M)$ is invariant under reparametrizations and under the group of isometries of $M$, and therefore yields an alternative to the elastic metric \eqref{eq:elasticmanifold}. The energy of a path of curves $[0,1]\ni t \mapsto c(t)$ for this metric is given by an expression similar to \eqref{eq:energy1},
\begin{equation}
\label{eq:energy2}
E(c) = \int_0^1\left(|\partial_tc(t,0)|^2+\int_I |\nabla_tq^\parallel(t,u)|^2\mathrm du\right)\mathrm dt.
\end{equation}
One finds that the conditions for such a curve to be a geodesic have been simplified with respect to those of the exact elastic metric framework written in Proposition \ref{prop:geod_equations_Riemannian}:
\begin{prop}[\cite{zhang2015video}]
A path of curves $[0,1]\ni t \mapsto c(t)$ is a geodesic minimizing the energy \eqref{eq:energy2} if and only if its TSRV representation $q^\parallel(t)=Q^\parallel(c(t))$ verifies the following equations:
\begin{align*}
\nabla_t\partial_tc(t,0) + \int_I \mathcal R(q^\parallel,\nabla_tq^\parallel)\partial_tc(t,u)\mathrm du  =& 0 , \quad \forall t\in[0,1], \\
\nabla_t^2 q^\parallel(t,u) =& 0 , \quad \forall (t,u)\in[0,1]\times I,
\end{align*}
where $\mathcal R$ denotes the curvature tensor of the base manifold $M$.
\end{prop}
In the context of finding the geodesic $c$ between two curves $c_1$ and $c_2$, the first equation describes the behavior of the baseline curve $t\mapsto c(t,0)$ linking the starting points $c_1(0)$ and $c_2(0)$, and the second equation expresses the fact that $q^\parallel=Q^\parallel(c)$ is covariant linear, i.e. $q^\parallel(t,u)$ can be obtained as a linear interpolation between the TSRV representations $q_1^\parallel(u)$ and $q_2^\parallel(u)$ of $c_1$ and $c_2$, parallel transported along the baseline curve to $c(t,0)$. The difficulty of implementing this method depends on the particular manifold $M$. For curves in the sphere $S^2$, the baseline curve linking the starting points is a circular arc, thus yielding simplifications with respect to the general geodesic shooting problem \cite{PhaseAmpSep}. The case of curves in the space of positive definite symmetric matrices is studied in \cite{zhang2018rate}.

\paragraph{Method 3.} A third possibility is to parallel transport the SRV representations of the curves to a particular reference point $p\in M$. This is the simplest method of all since the SRV representation of a curve is not only contained in a single linear space, but also this space is the same for all curves. The map of interest is then
\begin{equation*}
\begin{cases}
Q^{\parallel,p}:\operatorname{Imm}(I,M)\to C^\infty(I,T_pM)\\
Q^{\parallel,p}(c)=q^{\parallel,p},
\end{cases}
\end{equation*}
where $q^{\parallel,p}(u)$ is obtained by parallel translating the vector ${c'(u)}/{\sqrt{|c'(u)|}}$ along the shortest geodesic in $M$ from $c(u)$ to $p$. One then defines the distance between two curves $c_0$ and $c_1$ to be the $L^2$ distance between $q_0^{\parallel,p}=Q^{\parallel,p}(c_0)$ and $q_1^{\parallel,p}=Q^{\parallel,p}(c_1)$, i.e.
\begin{equation*}
    d(c_0,c_1)=\left(\int_I|q_0^{\parallel,p}(u)-q_1^{\parallel,p}(u)|^2\mathrm du\right)^{1/2}.
\end{equation*}
This distance function is invariant under reparametrizations of the curves, but it is not invariant under isometries of $M$. The main advantage of this method is computational speed. A disadvantage is that it depends heavily on the choice of the reference point $p$, and may induce serious distortions for curves that venture far away from  $p$. Finally, there can be problems with the definition of $Q^{\parallel,p}$ itself, since there can be more than one minimizing geodesic between $c(u)$ and $p$, and parallel translation along these different geodesics can yield different results.  In general, if all the curves being compared are not too far from the reference point $p$, this method can yield useful results at low computational cost; see \cite{SKKS2014} for applications to curves in $S^2$.

\section{Implementation}\label{sec:numerics}
In this section we will discuss the computation of the geodesic distance. We will first briefly address the case of parametrized curves. In the second part we will then describe the main difficulty in this context which is the minimization over reparametrizations in the group $\Diff(D)$. In particular we will describe several different approaches that have been developed to tackle this highly non-trivial task.

\subsection{The geodesic boundary value problem on parametrized curves}\label{min_parametrizedcurves}
 For open curves with values in Euclidean space, Lie groups or homogenous spaces and the SRV-metric, there exist analytic solution formulas for these operations and thus these computations become trivial. For most of the other situations discussed in this chapter, the absence of such formulas requires one to solve these problems using numerical optimization. Therefore, one first has to choose a discretization for all of the involved objects, i.e., one has to discretize the path of curves $c(t,u)$ for $t\in[0,1]$ and $u\in D$. A standard approach for this task consists of choosing B-splines in both time and space, i.e.
 \begin{equation}
c(t,u)=\sum_{i,j} c_{i,j} B_i(t) C_j(u)\;     
 \end{equation}
 where $B_i$ and $C_j$ are the chosen B-spline basis functions and where 
 $c_{i,j}$ for $i=0\ldots N_t$ and $j=0\ldots N_{u}$ are the coefficients.
 Note that this includes as a special case the discretization of regular curves as piecewise linear functions.
This procedure then reduces the
calculation of the geodesic distance~\eqref{distancefunction} to an unconstrained minimization problem of the discretized length functional, where the control points $c_{i,j}$ for $i=1\ldots N_t-1$ and $j=0\ldots N_{u}$
of the B-splines are the free variables. Here the control points of the boundary curves $c_{0j}$ and $c_{N_tj}$ are chosen as fixed parameters and are not changed in the optimization procedure. After this discretization step one can use standard methods of numerical optimization, such as the L-BFGS method, to approximate the solution of the finite-dimensional unconstrained minimization problem. For further information, in the notation of this chapter, we refer the reader to the article~\cite{BBHM2017}. See also~\cite{bauer2018relaxed,Vialard2014_preprint,Michor2006c}.

\subsection{Normalization by isometries}
The shape space $\mathcal S(D,M)$ in \eqref{eq:def_shape_space} involves quotienting out isometric transformations of $M$, in other words one has to technically minimize in \eqref{distancefunction_quotient} the elastic distance over $g \in \operatorname{Isom(M)}$. This is a finite-dimensional group which, for most manifolds $M$ encountered in practice, usually has a simple parametric representation.  

One common approach being used, although not rigorously equivalent to the optimization in \eqref{distancefunction_quotient}, is to pre-align the two shapes with respect to isometries of $M$ prior to estimating the elastic distance. When $M=\mathbb{R}^d$, this amounts to finding the optimal rotation and translation that best align them,  which is classically addressed by Procrustes analysis, cf. for example \cite{Drydmard16}. 

Alternatively, one can parametrize the group $\operatorname{Isom(M)}$ and perform the minimization over $g$ within the estimation of the distance itself, i.e. jointly with reparametrizations. For planar curves, this simply amounts to optimizing over a two-dimensional translation vector and the angle of rotation, which is the approach used, in particular, in \cite{BBHM2017,bauer2018relaxed}. Note that for general $\mathbb{R}^d$, a similar strategy is also possible by representing rotations as the exponential of antisymmetric matrices. In the case of manifold-valued curves however, normalizing with respect to isometries of $M$ may not always be relevant or can be harder to deal with in practice. This typically depends on the availability of convenient representations of the isometry group $\operatorname{Isom(M)}$, we refer the reader to \cite{su2018comparing} where some simple examples are considered.

\subsection{Minimization over the reparametrization group}
In addition to isometries of $M$, computation of distances and geodesics on the quotient space $\mathcal S(D,M)$ also requires to minimize the metric over reparametrizations in the group $\Diff(D)$, which is here infinite-dimensional. Several different approaches have been proposed to tackle this specific issue under various situations, which we review in the following paragraphs.

\subsubsection{Dynamic programming approach}\label{sec:dynamic}
A first method, which was proposed initially in \cite{trouve2000diffeomorphic,Mio2007}, is to convert this problem into a discrete optimization one. Considering piecewise linear (i.e. polygonal) curves, one may in turn choose to look for an optimal reparametrization of $\Diff(D)$ that is also piecewise linear. For curves in a Euclidean space and the SRV-metric, this is in part supported by the recent work of \cite{lahiri2015precise} where authors show that such optimal piecewise linear reparametrizations exist. In general, as piecewise linear functions are a dense set in the space of absolutely continuous functions, it is reasonable in practice to restrict the search to reparametrizations of this form. 

More specifically, assume that the two curves $c_0$ and $c_1$ are both piecewise linear. For simplicity, let's also assume that $D=[0,1]$ and that both curves are sampled uniformly on $D$, namely that $c_0$ and $c_1$ are linear on each of the subintervals $D_i = [t_i,t_{i+1}]$ for all $i=0,\ldots,N-1$ where $t_i=i/N$. One may then approximate positive diffeomorphisms in $\Diff(D)$ by piecewise linear homeomorphisms of $D$ with nodes in the set $\{0,t_1,t_2,\ldots,t_N\}$. Writing $J=\{t_0,t_1,t_2,\ldots,t_N\}$, we can equivalently consider all the polygonal paths defined on the grid $J\times J$ joining $(0,0)$ to $(1,1)$ and which are the graph of an increasing piecewise linear function with nodes in $J$. This set $\Gamma$ is now finite albeit containing a very large number of possible paths.

Nevertheless, an efficient way to determine an optimal discrete reparametrization is through dynamic programming. This is well-suited to situations where the energy to minimize can be written as an additive function over the different segments of the discrete path, which is made possible by the SRV transform in the case of elastic $G^{1,1/2}$ metrics (or more generally for the $G^{a,b}$-metric using the transforms of \cite{Michor2008a,KuNe2018,Bauer2014b}). We want to note here that this method is not well-suited to cases in which one does not has access to an explicitly computable distance function, such as for the higher-order elastic metrics.

Indeed, if $\gamma \in \Gamma$ is piecewise linear on the $K$ consecutive segments of vertices $(t_{i_0},t_{j_0})=(0,0)$, $(t_{i_1},t_{j_1}),\ldots,(t_{i_K},t_{j_K})=(1,1)$ with $t_{i_0}<t_{i_1}<\ldots<t_{i_K}$ and $t_{j_0}<t_{j_1}<\ldots<t_{j_K}$, then the discrete energy to be minimized is expressed as:
$$ E(\gamma) = \|Q(c_0) - Q(c_1\circ \gamma)\|_{L^2}^2 = \sum_{m=0}^{K-1} E(\gamma_{i_m,j_m}^{i_{m+1},j_{m+1}}) $$
where $E(\gamma_{i_m,j_m}^{i_{m+1},j_{m+1}})$ is the energy of the linear path from vertex  $(t_{i_m},t_{j_m})$ to $(t_{i_{m+1}},t_{j_{m+1}})$ and is given by
$$E(\gamma_{i_m,j_m}^{i_{m+1},j_{m+1}})=\frac{1}{N}\sum_{k=i_m}^{i_{m+1}-1} \left|Q(c_0)(t_k)-\sqrt{\frac{t_{j_{m+1}} - t_{j_m}}{t_{i_{m+1}} - t_{i_m}}}Q(c_1)(t_k) \right|^2.
$$ 
Now the generic dynamic programming method first computes the minimal energy among all paths in $\Gamma$ going from $(0,0)$ to any given vertex $(t_i,t_j)$, which we write $E^{i,j}$, through the following iterative procedure on $i$:
\vskip1ex
\noindent 1. Set $E^{(0,0)}=0$. \\
2. For a given $i\in\{1,\ldots,N\}$ and all $j\in \{1,\ldots,N\}$, compute $E^{(i,j)}$ and $P^{(i,j)}$ as:
\begin{equation}
\label{eq:dyn_prog_iteration}
E^{(i,j)} = \min_{(k,l) \in N_{ij}} E^{(k,l)} + E^{(i,j)}_{(k,l)}, \ \ P^{(i,j)} = \text{argmin}_{(k,l) \in N_{ij}} E^{(k,l)} + E^{(i,j)}_{(k,l)}
\end{equation}
where $E^{(i,j)}_{(k,l)}$ denotes in short the energy of the linear path from vertex $(t_k,t_l)$ to vertex $(t_i,t_j)$, and $N_{ij}$ is a set of admissible vertex indices connecting to $(i,j)$. 
\vskip1ex
\noindent At the end of this process, one obtains the minimal energy $E^{(N,N)}$. A corresponding optimal path $\gamma \in \Gamma$ can be simply recovered by backtracking from the final vertex $(1,1)$ to $(0,0)$, the index of the vertices in $\gamma$ being specifically $(i_q,j_q)=(N,N)$, $(i_{q-1},j_{q-1}) = P^{(i_q,j_q)},\ldots,(i_1,j_1)=P^{(i_{2},j_{2})}$ and $(i_0,j_0)=P^{(i_1,j_1)}=(0,0)$. 

The choice of search neighborhood $N_{ij}$ in the above procedure has a critical impact on the resulting complexity. To find the true minimum over all possible paths in $\Gamma$, one should technically take in \eqref{eq:dyn_prog_iteration}, $N_{ij} = \{(k,l):0\leq k \leq i-1, \ 0\leq k \leq j-1 \}$ for any $1\leq i,j \leq N-1$. This would result however in a high numerical cost of the order $O(N^4)$. It can be significantly reduced by restricting $N_{ij}$ to a smaller set of admissible neighboring vertices. For instance, authors in \cite{Mio2007} propose to limit the search to a small square of size $3\times 3$ with upper right vertex $(i-1,j-1)$.
While this constrains the possible minimal and maximal slope of the estimated $\gamma$, it is generally sufficient in most cases and reduces the numerical complexity to $O(N^2)$, making the whole approach efficient in practice. Note that alternative dynamic programming algorithms have been investigated more recently, in particular in the work of \cite{Bernal2016} which makes use of adaptive strips neighborhoods to further reduce the complexity to $O(N)$.   

\subsubsection{Discretizing the diffeomorphism group and using gradient based methods}\label{sec:discretizediff}
A second method, which has been proposed in the context of the SRV-metric in \cite{huang2016riemannian,huang2014riemannian} and for higher-order Sobolev metrics in \cite{BBHM2017}, is also based on a direct discretization of the diffeomorphism group and the space of curves. However, in contrast with the previous section where diffeomorphisms of $D$ were discretized as piecewise linear functions, this method offers more flexibility. For example one could choose -- similarly to Section~\ref{min_parametrizedcurves} -- B-spline representations of reparametrizations. Considering the distance function~\eqref{dist_shape} on the space of unparametrized curves in this discretization leads again to a finite-dimensional minimization problem, which can be tackled by standard methods. 

In the case when one has no access to an explicit formula for the geodesic distance -- such as for higher-order Sobolev metrics -- it is computationally efficient to view this problem as a joint minimization problem over the (discretized) path of curves
$$c(t,u)=\sum_{i,j} c_{ij}B_i(t)C_j(u)$$ and the reparametrization function $$\gamma(u)=\sum_k \gamma_k D_k(u).$$ Here $B_i, C_j$ and $D_k$ are the chosen basis functions for the discretization of the path of curves and the reparametrization function respectively. 
One difficulty in this context is that the composition of the (discretized) target $c(1,u)$ and the (discretized) reparametrization function $\gamma(u)$ typically leaves the chosen discretization space. Thus one has to consider the corresponding projection operator that projects this reparametrized curve back to the discretization space. This procedure can lead to numerical phenomenona such as loss of features in the target curve. For more details we refer to the presentation in \cite{BBHM2017}.

\subsubsection{Iterative ``horizontalization" method}\label{sec:horizontalizations}

\begin{figure}
    \centering
    \includegraphics[width=18em]{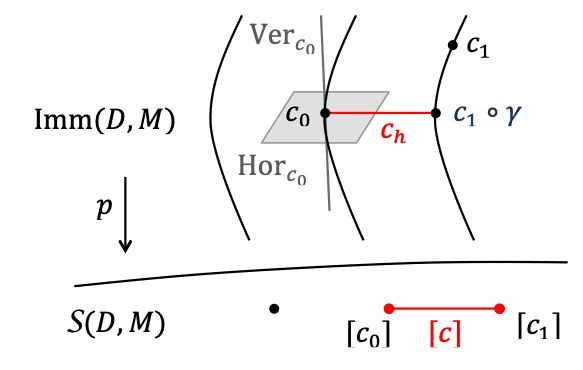}
    \caption{Principal bundle structure formed by the space of curves and their shapes. The horizontal geodesic $c_h$ between $c_0$ and the optimally matched $c_1\circ\gamma$ projects to a geodesic $[c]=p(c_h)$ between the corresponding shapes.}
    \label{fig:principal_bundle}
\end{figure}

Another possibility is to exploit the principal bundle structure formed by the space of parameterized curves and their shapes. The fibers of this bundle are the sets of all the curves that are identical modulo reparametrization, i.e. that project onto the same shape (Figure \ref{fig:principal_bundle}). Any tangent vector $h\in T_c\Imm(D,M)$ can be decomposed as the sum of a vertical part $h^{ver}\in \text{Ver}_c$ tangent to the fiber, which has an action of reparameterizing the curve without changing its shape, and a horizontal part $h^{hor}\in\text{Hor}_c = \left(\text{Ver}_c \right)^{\perp_G}$, $G$-orthogonal to the fiber. While the horizontal subspace depends on the choice of the reparametrization invariant metric $G$, the vertical subspace is always the same:
\begin{equation*}
\text{Ver}_c = \ker dp(c) = \left\{ m v := m c'/|c'| : \,\, m \in C^\infty([0,1], \mathbb R),\, m(0)=m(1)=0 \right\}.
\end{equation*}
Paths of curves with horizontal velocity vectors are called horizontal, and horizontal geodesics for $G$ project onto geodesics of the shape space for the Riemannian metric induced by the Riemannian submersion $p:\Imm([0,1],M)\rightarrow \mathcal S([0,1],M)$, see e.g. \cite[Section 26.12]{michor2008topics}. A natural way to solve the boundary value problem in the shape space is by fixing the parameterization $c_0$ of one of the curves and computing the horizontal geodesic linking $c_0$ to the closest reparametrization $c_1\circ\gamma$ of the second curve $c_1$, by iterative "horizontalizations" of geodesics. The idea is to decompose any path of curves $t\mapsto c(t)\in\Imm(D,M)$ as
\begin{equation}
\label{def}
c(t,u) = c^{hor}(t,\gamma(t,u)) \quad \forall (t,u)\in[0,1]\times D,
\end{equation}
where $t\mapsto c^{hor}(t)$ is a horizontal path and is reparameterized by a path of diffeomorphisms $t\mapsto \gamma(t) \in \text{Diff}^+(D)$.
Differentiating with respect to $u$ and $t$ and taking the squared norm with respect to $G$ yields 
\begin{align*}
    |\partial_uc|^2 &= |\partial_u\gamma|^2 |\partial_uc^{hor}\circ\gamma|^2,\\
    |\partial_tc|^2&=|\partial_tc^{hor}\circ\gamma|^2+|\partial_t\gamma|^2|\partial_uc^{hor}\circ\gamma|^2,
\end{align*}
where in the second expression we have used the fact that $\partial_tc^{hor}\circ\gamma$ is horizontal by definition of $c^{hor}$, and $\partial_uc^{hor}$ is vertical as we can see from the first expression. From this, we immediately see that if the metric $G$ is reparameterization invariant, taking the horizontal part of a path decreases its length
\begin{equation*} L_G(c^{hor}) \leq L_G(c).
\end{equation*}
Therefore, by taking the horizontal part of the geodesic linking two curves $c_0$ and $c_1$, we obtain a shorter, horizontal path linking $c_0$ to the fiber of $c_1$, which gives a closer (in terms of $G$) representative $\tilde c_1 = c_1\circ \gamma(1)$ of the target curve. However it is no longer a geodesic path. By computing the geodesic between $c_0$ and this new representative $\tilde c_1$, we are guaranteed to reduce once more the distance to the fiber. The optimal matching algorithm simply iterates these two steps, and converges to a horizontal geodesic. At each step, the horizontal part of the geodesic can be computed using the following result.
\begin{prop}[\cite{lebrigant2019discrete}]
The path of diffeomorphisms $t\mapsto\gamma(t)\in\Diff(D)$ that transforms a path $t\mapsto c(t)\in\Imm(D,M)$ into a horizontal path is solution of the PDE
\begin{equation}
\label{phicont}
\partial_t\gamma(t,u) = \frac{m(t,u)}{|\partial_uc(t,u)|} \partial_u\gamma(t,u),
\end{equation}
with initial condition $\gamma(0)=\text{Id}$, and where $m(t,u):=|\partial_tc^{ver}(t,u)|$. %
\end{prop}
This method can be applied as long as the horizontal part of a tangent vector (or equivalently, the norm of the vertical component $m$) can be computed. For the class of $G^{a,b}$-elastic metrics, and for the SRV-metric in particular, $m$ can be found by solving an ODE, see \cite{lebrigant2019discrete}. An example of geodesic between curves in the hyperbolic plane estimated with this approach is shown in Figure \ref{fig:example_curve_geodesics} (right).

\begin{figure}
    \centering
    \begin{tabular}{cccc}
    \includegraphics[width=3.3cm]{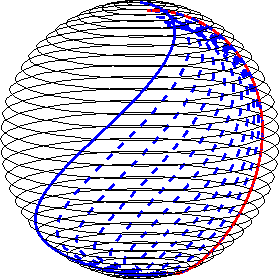} & \includegraphics[width=3.3cm]{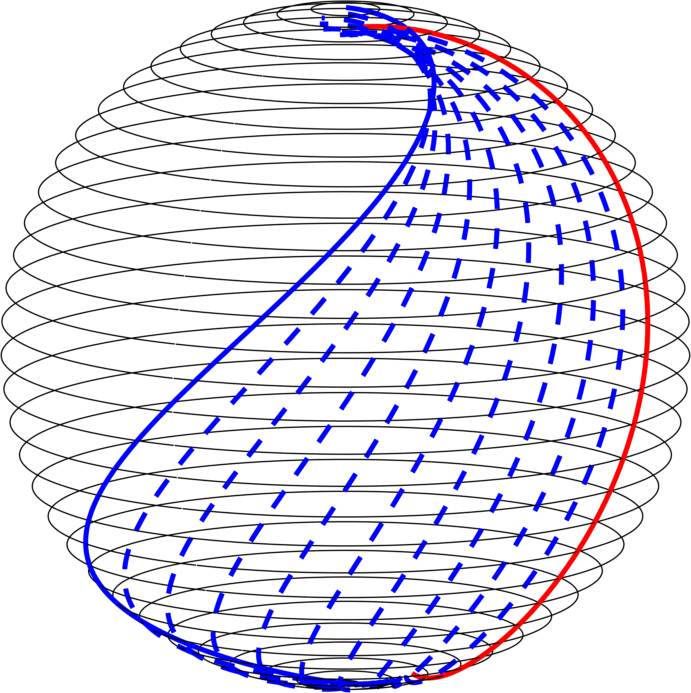} & \includegraphics[width=3.3cm]{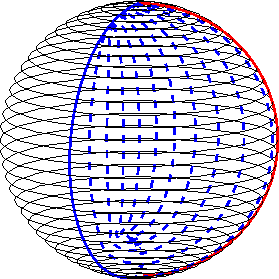}&
    \includegraphics[width=3.3cm]{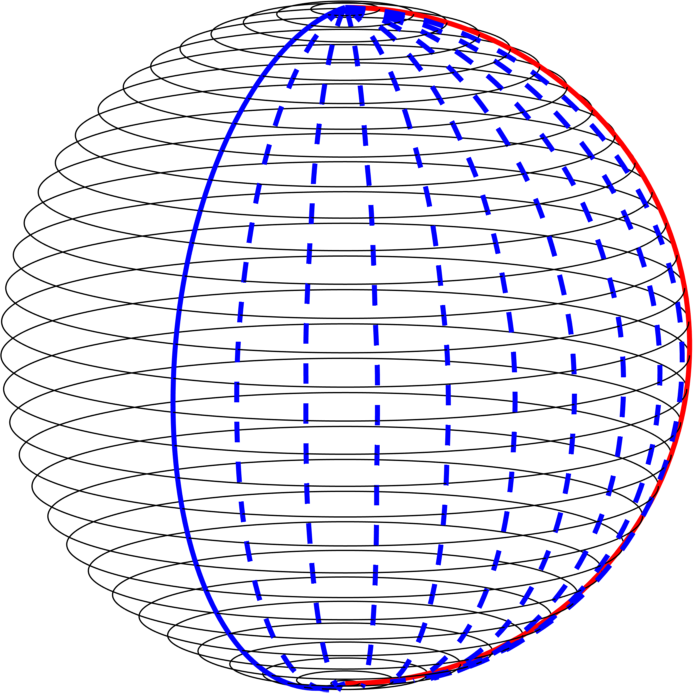}
    \\
    $\text{dist}=1.32$ & $\text{dist}=2.03$ &  $\text{dist}=1.92$&$\text{dist}=2.47$
    \end{tabular}
    \caption{Numerical comparison of the distance and geodesics between 3D curves lying on the unit sphere modulo reparametrizations: first and third picture for the SRVF metric in the Euclidean space $\mathbb{R}^3$ (computed with the relaxed algorithm of \cite{Bauer_GSI2019}), second and forth picture for the SRVF distance on $S^2$ (estimated with the method of \cite{su2018comparing}). Observe that the geodesics calculated in Euclidean space do not stay on the sphere and thus result in a lower SRVF distance.}
    \label{fig:example_curve_geodesics_3D}
\end{figure}

\subsubsection{Relaxation of the exact matching problem}
\label{ssec:relaxed_approach}
A last possible approach to deal with reparametrization invariance in the computation of geodesics and distances on the quotient space (without directly optimizing over reparametrizations) is to introduce a relaxation term for the end time constraint providing a measure of discrepancy up to reparametrization to the target curve $c_1$. This is inspired by similar methods used earlier on in diffeomorphic registration frameworks, see e.g. \cite{Glaunes2008,Durrleman2010,Charon2013,Roussillon2016,Charon2017} among other references. But it can also be applied in the context of elastic metric matching, as recent works such as \cite{bauer2018relaxed,Bauer_GSI2019,sukurdeep2019inexact} have shown. In this section, we will assume that curves are immersed in the Euclidean space $\mathbb{R}^d$.

Going back to the original formulation of the geodesic distance given by \eqref{distancefunction} and \eqref{distancefunction_quotient}, the idea is to start by replacing the end time boundary constraint that $c(1)=c_1 \circ \gamma$ for some $\gamma \in \Diff(D)$ using a surrogate fidelity (or discrepancy) term $\tilde{d}(c(1),c_1)$. Assuming that $\tilde{d}(c(1),c_1)$ is invariant to the parametrization of \textit{both} $c(1)$ and $c_1$, i.e. that $\tilde{d}$ defines a distance on the quotient space, one gets the equivalence  between the above boundary condition and $\tilde{d}(c(1),c_1) = 0$. Then we may choose to relax the constraint and consider the alternative variational problem: 
\begin{equation}
\label{eq:geod_relaxed}
  \operatorname{inf} \int_0^1 G_c(\partial_t c,\partial_t c) dt + \lambda \tilde{d}(c(1),g\circ c_1)^2 
\end{equation}
over all paths $c: [0,1]\to\Imm(D,\mathbb{R}^d)$ such that $c(0)=c_0$. Note that minimization over $\gamma \in \Diff(D)$ is no longer needed here and a minimizing path $c$ of \eqref{eq:geod_relaxed} is by construction a geodesic between $c_0$ and $c(1) \approx c_1$ in the quotient space $\mathcal S(D,\mathbb{R}^d)$. In the above, $\lambda>0$ denotes a fixed weighting coefficient between the two terms which controls the accuracy of the matching to the target $c_1$. Other strategies such as augmented Lagrangian methods can also be used to adapt the choice of this parameter in order to reach a prescribed matching accuracy, cf. \cite{bauer2018relaxed}.

\begin{rem}
In the specific case of the SRV-metric of Section \ref{ssec:SRV_euclidean}, the variational problem \eqref{eq:geod_relaxed} can be even further simplified to a minimization problem over the end curve $c^1=c(1) \in \Imm(D,\mathbb{R}^d)$ instead of a full curve path. Indeed, using the properties of the SRV transform, it is easy to see that the problem can be equivalently rewritten as:
\begin{equation*}
   \operatorname{inf}_{c^1} \|Q(c^1) - Q(c_0)\|_{L^2}^2 + \lambda \tilde{d}(c^1,g\circ c_1)^2.   
\end{equation*}
and leads, after discretization, to a simple minimization problem over the vertices of the deformed curve. This formulation is for instance implemented in \cite{Bauer_GSI2019}. Note that this principle also applies to other simplifying transforms associated to different choices of elastic parameters as proposed and implemented in \cite{sukurdeep2019inexact}.
\end{rem}

This entire approach relies on the discrepancy distance $\tilde{d}$ which, in particular, needs to be itself independent of curve parametrization. This may sound redundant as this is also the purpose of the quotient metric construction we have been discussing all along in this chapter. Yet one can construct discrepancy metrics that are both simple and easy to compute in practice, i.e. that do not require solving an extra optimization problem. Even though these discrepancy distances do not fit within the Riemannian metric setting that we are ultimately interested in, they remain ideally suited as auxiliary terms within the elastic matching problem. While different constructions are possible, the key strategy developed in the aforementioned references consists in embedding any unparametrized curve into a certain measure space and thereby recover explicit distances derived from kernel metrics on this measure space. We will however not elaborate on the actual construction of such embeddings and metrics; the interested reader may refer to the recent survey of \cite{Charon2020}.

Unlike the methods discussed in the previous sections, this relaxed approach does not necessarily compute the exact distance between the two curves. Yet it can prove particularly useful in situations where one or both curves are corrupted by noise or small topological perturbations that may otherwise considerably affect the estimated value of the distance. In addition to the example of Figure \ref{fig:example_curve_geodesics} (left), we show in Figure \ref{fig:example_curve_geodesics_3D} additional geodesics for the SRVF metric between curves of $\mathbb{R}^3$ (lying on the unit sphere) estimated by this approach, which we compare to the geodesics for the SRVF metric on the homogeneous space $S^2$.     

\subsection{Open source implementations}
Several of the methods and algorithms described above are available in open source software packages. Here is a (non-exhaustive) list of some of these:
\begin{itemize}
    \item {\bf Second order elastic metrics for curves in $\mathbb R^d$:} Implementation of a four-parameter family of metric (including in particular the family of $G^{a,b}$-metric) is available at:
        \begin{center}
    \url{https://github.com/h2metrics/h2metrics}
    \end{center}
    Both the inexact matching approach of Section~\ref{ssec:relaxed_approach} and the gradient based approach of Section~\ref{sec:discretizediff} are implemented.
    \item {\bf SRV framework for curves in $\mathbb R^d$}: several different implementations for this classical method exist. 
    This includes in particular the R-package by J. Tucker
        \begin{center}
    \url{https://cran.r-project.org/web/packages/fdasrvf/}
    \end{center}
    and 
    the Matlab implementation of M. Bruveris as available on github:
    \begin{center}
        \url{https://github.com/martinsbruveris/libsrvf}
    \end{center} 
    In the second one, both the dynamic programming approach of Section~\ref{sec:dynamic} and the explicit solution formula discussed in Remark~\ref{rem:explicit_optimalreparams} are implemented. 
    
    \item {\bf SRV metric for curves in homogenous spaces and Lie groups:} Code for several choices for the target space $M$ can be found at:
            \begin{center}
    \url{https://github.com/zhesu1/SRVFhomogeneous}
    \end{center}
    Optimal reparametrizations are estimated using the dynamic programming approach of Section~\ref{sec:dynamic}.
\end{itemize}

\section{Conclusion}
In this chapter, we reviewed the current state-of-the-art of curve comparison through intrinsic quotient Riemannian metrics for Euclidean as well as non-Euclidean curves. We discussed the theoretical framework, in particular the questions of non-degeneracy of Sobolev metrics and geodesic completeness of the corresponding infinite-dimensional manifolds before analyzing more specifically the case of the SRV-metric for which the variational expression of the distance considerably simplifies. We also discussed several numerical approaches that have been proposed for the computation of such metrics in the different settings and for which several open source implementations are available.

There are many directions in which this framework can be extended. One is the construction and computation of corresponding intrinsic metrics between surfaces modulo reparametrizations. Due to their significantly more complex structure than curves, this is a subject of ongoing and active investigations both from the mathematical and numerical side: we refer interested readers e.g. to \cite{jermyn2017elastic,Kurtek2011,su2019shape_surfaces,tumpach2015gauge,kilian2007geometric}.

Going back to curves, as noted in Remark \ref{rem:simplifying_transforms}, there have been several extensions and variations of the SRV framework which introduced simplifying transforms for other first-order metrics than the specific one considered in Section \ref{ssec:SRV_general}. We finally mention the recent work of \cite{Younes2018_hybrid} which explored the possibility to combine intrinsic Sobolev metrics with extrinsic diffeomorphism-based metrics within a hybrid framework.       

\section*{Acknowledgements}
M. Bauer was partially supported by NSF-grant 1912037 (collaborative research in connection with NSF-grant 1912030) and NSF-grant 
1953244 (collaborative research in connection with NSF-grant 1953267).
N. Charon was partially supported by NSF-grant 1945224 and NSF-grant 1953267 (collaborative research in connection with NSF-grant 1953244). Eric Klassen gratefully acknowledges the support of the Simons Foundation-grant 317865.

\end{document}